\documentclass[reqno]{amsart}
\usepackage{epsfig, amsmath, amsfonts, amsthm, amssymb, amscd, eufrak, latexsym}
\usepackage[all]{xy}
\begin{document}
\newcommand{\motclef}[1]{{\em #1}}
\newtheorem{theo}{Theorem}[section]
\newtheorem{prop}{Proposition}[section]
\newtheorem{definition}{Definition}[section]
\newtheorem{expl}{Example}[section]
\newtheorem{remark}{remark}[section]
\newtheorem{lem}{lemma}[section]
\newtheorem{cor}{Corollary}[section]
\newtheorem{conj}{Conjecture}[section]
%%%%%%%%%%%%%%%%%%%%%%%%%%%%%%%%%%%
%%%%%%%%%%%%%%%%%%%%%%%%%%%%%%%%%%%%
%
%
\thanks{ }
\title[ON REPRESENTATIONS OF BOL ALGEBRAS]{ON REPRESENTATIONS OF BOL ALGEBRAS}
\author[N. Ndoune]{Ndoune Ndoune}
\address{Ndoune Ndoune \newline
D\'epartement de Math\'ematiques,
\newline
Universit\'e de Sherbrooke,
\newline
Sherbrooke, (Qu\'ebec),
\newline
J1K 2R1, CANADA.
\newline {}}
\email{ndoune.ndoune@usherbrooke.ca}
\author[T. Bouetou]{Thomas Bouetou Bouetou}
\address{Thomas Bouetou Bouetou\newline
Department of Mathematics and Computer Science,
\newline
Ecole Nationale
Sup\'erieure Polytechnique de
Yaound\'{e},
\newline
 and, Centre d'Excellence Africain en Technologies de l'Information et de la Communication (CETIC),\newline
 the Abdus Salam International Centre for Theoretical Physics (ICTP)
\newline
P.O.Box: 8390 Yaounde CAMEROUN},
\email{tbouetou@ictp.it}
\maketitle

\begin{abstract}
In this paper, we introduce the notion of representation of Bol algebra.
We prove an analogue of the classical Engel's theorem and the extension
of Ado-Iwasawa theorem for Bol Algebras. We study the category of
representations of Bol algebras and show that it is a tensor category.
In the case of regular representations of Bol algebras, we give a
complete classification of them for all two-dimensional Bol algebras.
\end{abstract}

\section{Introduction}

It is well known that the algebraic systems which characterize
locally a totally geodesic subspace is a Lie triple system see
(\cite{Car}, \cite{Nom}, \cite{Yam}). A Bol algebra is
realized by equipping Lie triple System with an additional binary
skew operation which satisfies a pseudo-differentiation property
(\cite{kuza}, \cite{pi}). A morphism of Bol algebras is a linear map
which preserves the ternary and the binary
operations. More generally, the algebraic structures
which characterize locally Bol loops are Bol algebras \cite{sa1}.
Until now, the representations of these algebras have not been studied.
Since the representations of Lie algebras and Lie groups have natural
connection with particulars physics, we claim that the representations
of Bol algebras should lead with the physical applications. More precisely,
in physics the representations of Bol algebras will be useful for
the description of invariant properties of physical systems.
and the concomitant conservation laws as a result.
In \cite{mop} it is shown that, Malcev algebras and
 Lie triple systems are particular subclasses of Bol algebras.
The representations of Malcev algebras can be found in
\cite{ku}, and those of Lie triple systems was constructed
in \cite{hob} and \cite{wom}. Now, there already exists some representations
of other classes of non-associative algebras; the case of
alternative algebras was constructed by Schafer \cite{Shaf},
the one of Leibniz algebras by  Kolesnikov \cite{Kol} and for
Jordan superalgebras, the representations was given by  Consuelo and
 Zelmanov see \cite{maze}.\\
Let $\mathfrak{B}$ be a Bol algebra over a field $K$ of characteristic zero,
a representation
of Bol algebra $\mathfrak{B}$ on a $K$-vector space $V$ is a triplet of maps
$(\rho,\delta,\Delta)$ which respect some conditions which will be given
later in
the paper.\\ Our first main result is the following.

\begin{theo}
Let $\mathfrak{B}$ be a finite dimensional Bol algebra over a field $K$
and $\mathcal{R}$ consist of nilpotent representations of Bol algebra
$\mathfrak{B}$ in a finite dimensional space $V$. Then there exists a
vector $v\in V$, $v\neq 0$ such that $(\rho,\delta,\Delta)(v)=0$ for
all $(\rho,\delta,\Delta)\in
\mathcal{R}.$

\end{theo}

We agree that the image of any vector $v$ of $V$ by the operator $(\rho,
\delta,\Delta)$ is given by $(\rho,\delta,\Delta)(v)=(\rho(v_1),
\delta(v_2),\Delta(v_3))$, where $v=(v_1,v_2,v_3)\in \mathfrak{B}^{3}$.\\
We define also the regular representations and the
adjoint representations of Bol algebras. As
an easy consequence, we show that if any representation of Bol algebra is nilpotent,
 then its adjoint representation is also nilpotent.\\
 We are also interested by the question of the extension theorem of Ado-Iwasawa
  for Bol algebras. In \cite{pi} P\'erez-Izquierdo established the existence of a
  Poincar\'e-Birkhoff-Witt type basis for a universal envelope of Bol algebra.
The above result allows us to interest ourselves to an extension of Ado-Iwasawa theorem for Bol algebra.
let $A$ be an alternative algebra, the  the generalized
right alternative nucleus is the algebra $RN_{alt}(A)$ defined by
$RN_{alt}(A)=\{a\in A/ (x,a,y)=-(x,y,a)\}$ .
We then give our second theorem.

\begin{theo}
Let $\mathfrak{B}$ be a finite-dimensional right Bol algebra over a field of characteristic
different to 2 and 3. Then there exists a unital finite-dimensional algebra $A$ and a
monomorphism of Bol algebras $\mathfrak{B}\longrightarrow RN_{alt}(A).$
\end{theo}
The analogue of our second result above was established for Malcev algebras in \cite{izsh}.
The collection of all representations of Bol algebra and the morphisms between them form a category, named
 the category of representations of Bol algebras $Rep(\mathfrak{B}).$ One can view a representation of Bol
 algebra as a $\mathfrak{B}$-module analogously as in \cite{maze} in the case of Jordan superalgebras.
 One can understand also the representations of Bol algebras in term of matrices with sweet properties.
 The investigation between
 the category $Rep(\mathfrak{B})$ and the category of left $U(\mathfrak{B})$-modules, where $U(\mathfrak{B})$
 is the universal enveloping algebra of $\mathfrak{B}$, endowed with its bialgebra structure, leads
 us to our third main theorem.

\begin{theo}
The category of representations of Bol algebra $Rep(\mathfrak{B})$ is equivalent to the category of
representations of its universal enveloping algebra $Rep(U(\mathfrak{B})).$
\end{theo}

The paper is organized as follows: We introduce in section 2 the notion of representations of Bol algebra.
In section 3 we establish the Engel's theorem for Bol algebras.
In section 4 an extension of Ado-Iwasawa theorem to Bol algebras is proved. Finally in section 5, we present
the category
of representations of Bol algebras and show that it is equivalent to the category of left modules under its
universal enveloping algebra. As immediate consequence, we show the category $Rep(\mathfrak{B})$ is a tensor
 category. We end the section by given a complete classification of
regular representations of two-dimensional Bol algebras.

\section{Bol algebras and their representations}

Bol algebras were introduced in differential geometry to study smooth Bol loops (\cite{sa1},\cite{sa2},\cite{sami}).
 A $right$ $loop$ is a set $\mathcal{Q}$, together with a binary operation $(a,b)\mapsto a\cdot b$, such that for
 any $b$ in $\mathcal{Q}$, the right multiplication operator $R_{b}:x\mapsto x\cdot b$ is bijective, and there exists
 an element $\varepsilon \in \mathcal{Q}$, such that $\varepsilon\cdot b=b$ for any $b$ in $\mathcal{Q}$. The dual
 definition gives rise to a left Bol loop. In case that $\langle \mathcal{Q},\cdot,\varepsilon\rangle$ is both left
and right loop then it is called a loop with identity element $\varepsilon.$ \newline
A $right$ $smooth$ $loop$ $\mathcal{M}$ is a right loop equipped with a structure of smooth manifold, that is the map
$(a,b)\mapsto a\cdot b$ and $R_{b}^{-1}$ are smooth, see (\cite{sa2},\cite{sami}.) Since groups are particular loops,
so the Lie groups are particular cases of smooth loops. In scientific literature, many classes of loops are known:
homogeneous loops, Moufang loops, Bol loops, Kikkawa loops among others.\newline
A $right$ $Bol$ $loop$ $\langle \mathcal{Q},\cdot,\varepsilon\rangle$ is a right loop that satisfies the right Bol
identity $$x\cdot ((a\cdot y)\cdot a)=((x\cdot a)\cdot y)\cdot a$$ for all $a,x,y$ in $\mathcal{Q}$. Similarly, a
$left$ $Bol$ $loop$ satisfies the identity $a\cdot (x\cdot (a\cdot y))=(a\cdot (x\cdot a))\cdot y$.

As in the case
of Lie groups where the tangent space at each point is equipped with Lie algebra structure, the tangent space at each
point of Bol loop is equipped with the structure of Bol algebra.

\begin{definition}
A vector space $\mathfrak{B}$ over a field $K$
is called Bol algebra if it is equipped with a trilinear operation $[-; -,-]$ and
a skew-symmetric operation $x \cdot y$
satisfying the following identities: \begin{enumerate} \item[(i)] $[x;x,y]=0$ \item[(ii)]
$[x;y,z]+[z;x,y]+[y;z,x]=0.$ \item[(iii)]
$[[x; y,z];\alpha,\beta]=[[x;\alpha,\beta];y,z]+[x;[y;\alpha,
\beta],z]+[x;y ,[z;\alpha,\beta]]$ \item[(iv)]
$[x \cdot y;\alpha,\beta]=[x;\alpha,\beta]\cdot y+
x \cdot [y;\alpha,\beta]+ [\alpha \cdot \beta; x,y]+[x \cdot
y]\cdot [\alpha \cdot \beta]$
\end{enumerate}
for all $x$,$y$,$z$,$\alpha$ and $\beta$ in
$\mathfrak{B}$.\end{definition}

In other words, a Bol algebra is a
Lie triple system $(\mathfrak{B},[-;-,-])$ with an additional
bilinear skew-symmetric operation  $x \cdot y$ such that, the
derivation $D_{\alpha,\beta}:x\longrightarrow
[x;\alpha,\beta]$ on a ternary operation is a pseudo-differentiation
with components $\alpha,\beta$ on a binary operation, that is; for
all $x$,$y$ and $z$ in $\mathfrak{B}$, we have
$$D_{\alpha,\beta}(x \cdot y)=\big(D_{\alpha,\beta}(x) \big)\cdot y+
x \cdot \big(D_{\alpha,\beta}(y) \big)+ [\alpha \cdot
\beta; x,y]+(x \cdot y)\cdot (\alpha \cdot \beta).$$

$D_{\alpha,\beta}$ is a differentiation on ternary
operation $[-;-,-]$ that is;
$$D_{\alpha,\beta}[x; y,w]=[D_{\alpha,
\beta}(x);y,w]+[x;D_{\alpha, \beta}(y),w]+[x;y
,D_{\alpha,\beta}(w)].$$
In fact, the Bol algebra
defined above is called right Bol algebra see \cite{pi}.
In particular, any Lie triple system may be regarded as Bol algebra
with the trivial multiplication $x\cdot y=0$, for all
$x,y \in \mathfrak{B}.$\\

Bol algebras can be realized as the tangent
algebras of  Bol loops with the right  Bol identity, and they
allow embedding in Lie algebras \cite{sa1}, \cite{sa2}.

\begin{definition}

A linear map $\varphi: \mathfrak{B_{1}}\rightarrow
 \mathfrak{B_{2}}$ between two Bol
algebras is called morphism of Bol algebras if it is preserve
the ternary and the binary
operations.
\end{definition}

The subspace $S$ of Bol algebra $\mathfrak{B}$ is a sub-Bol algebra
if the inclusion $j:S\hookrightarrow \mathfrak{B}$ is a morphism
of Bol algebras.

\begin{definition}
Let $(\mathfrak{B},[-;-,-],\cdot)$ be a Bol algebra
over a field $K$, a pseudo-differentiation
is a linear map $P:\mathfrak{B}\longrightarrow \mathfrak{B}$ for
which, there exists $z\in \mathfrak{B}$ (a companion of $D$) such
that $P(x\cdot y)=P(x)\cdot y+x\cdot P(y)+[z;x,y]+(x\cdot y)\cdot
z;$ the companion is not necessarily unique.
\end{definition}

The set of all companions of $D$ is denoted $Com(D).$
The map $D_{\alpha,\beta}:x\longrightarrow
[x;\alpha,\beta]$ is a pseudo-differentiation with companion
$\alpha\cdot \beta$, called inner pseudo-differentiation of $\mathfrak{B}$.
The pseudo-differentiations of $\mathfrak{B}$ form a Lie algebra,
denoted by pder$\mathfrak{B}$ under the natural product  $[P,P']=PP'-P'P$.
The algebra ipder$\mathfrak{B}$ generate by $\{D_{a,b}/ a,b\in \mathfrak{B}\}$
is a Lie subalgebra of pder$\mathfrak{B}$, called the Lie algebra of inner
pseudo-differentiations of $\mathfrak{B}$.
The enlarged algebra of  pseudo-differentiations of $\mathfrak{B}$
is defined as $Pder\mathfrak{B}=\{(D,z),D\in pder\mathfrak{B},z\in
Com(D) \}$\\
and the enlarged algebra of inner
pseudo-differentiation is defined as
$Ipder\mathfrak{B}=\{(D,z),D\in ipder\mathfrak{B},z\in Com(D) \}$.\\
It is showed in \cite{kuza} and \cite{pi} that, those algebras
defined below are the Lie algebras with the brackets $[P,P']=PP'-P'P$.\\
The direct sum $L=\mathfrak{B}\oplus Ipder\mathfrak{B}$ is a Lie algebra with
the operation $[x,y]=D_{x,y}$, $[x,D_{a,b}]=D_{a,b}(x)$, for all $x,y,a,b$ in
$\mathfrak{B}.$ The Lie algebra $(L,[,])$ is called the standard enveloping
Lie algebra of Bol algebra
$\mathfrak{B}.$
\\
The map $\delta_a:x\longmapsto x\cdot a$ is
a linear map of $\mathfrak{B}$. We denote by $\overline{\mathfrak{B}}$ the Lie
algebra generate by $\{\delta_a, a\in \mathfrak{B}\}$  with brackets
$[\delta_a,\delta_b]=\delta_a\delta_b-\delta_b\delta_a.$ We get an other Lie
 algebra $\overline{L}=\overline{\mathfrak{B}}\oplus Ipder\mathfrak{B}$
 which is a subalgebra of the Lie algebra generated by linear maps
 of $\mathfrak{B}$.\\
If the subspace $\mathcal{I}$
of $\mathfrak{B}$ satisfies the stronger condition $\mathcal{I}\cdot\mathfrak{B}+( \mathcal{I};
 \mathfrak{B},\mathfrak{B})\subset \mathcal{I}$, then  $\mathcal{I}$ is an ideal of $\mathfrak{B}.$
  An ideal $\mathcal{I}$ of $\mathfrak{B}$ automatically satisfies $(\mathfrak{B};\mathcal{I},
  \mathfrak{B})\subset \mathcal{I}$ and $(\mathfrak{B};
 \mathfrak{B},\mathcal{I})\subset \mathcal{I}.$\\
 For more understanding of Bol algebras and Bol loops, it is important to investigate about
 their representations.
 We defined a representation of
  Bol algebra as follows.

\begin{definition}
If $\mathfrak{B}$ is a Bol algebra over a field $K$ and $V$ a vector field
 over $K$, the pair $(\rho,\delta)$ with the skew-symmetric bilinear map
 $\rho:\mathfrak{B}^{2}\longrightarrow EndV$ and the linear map
 $\delta:\mathfrak{B}\longrightarrow EndV$ is said to be a representation
 of Bol algebra $\mathfrak{B}$ in $V$ if there exists a bilinear operation
$\Delta:\mathfrak{B}^{2}\longrightarrow EndV$ such that the following statements
are satisfied:

 \begin{enumerate} \item[(R1)] $\rho(u,v)=\Delta(u,v)-\Delta(v,u)$ \item[(R2)]
$[\rho(a,b),\rho(u,v)]=\rho([a,u,v],b)+\rho(a,[b,u,v])$ \item[(R3)]
$[\rho(u,v),\delta(a)]=\delta([a,u,v])+\Delta(u\cdot v,a)+\delta(u\cdot v)\delta(a)$

\end{enumerate}
for all $x,y,a,b$ in $\mathfrak{B}$.

\end{definition}
The operation $\Delta$ is called a companion of the representation
$(\rho,\delta)$ of the Bol algebra $\mathfrak{B}$.\\
 In this case we
 can
denoted by $(\rho,\delta,\Delta,V)$ or simply $(\rho,\delta,\Delta)$, the representation $(\rho,\delta,V)$ with
companion $\Delta.$ Following the approach of Consuelo and Zelmanov
for the representations of Jordan Superalgebras in \cite{maze},
it is equivalent to say that the vector space $V$ is a Bol module
($\mathfrak{B}$-module) i.e., $E_{V}=\mathfrak{B}\oplus V$ possesses
the structure of Bol algebra such that:

\begin{enumerate} \item
[(a)] $\mathfrak{B}$ is a sub-Bol
algebra of $E_{V}$, \item
[(b)] $V$ is an ideal of Bol algebra $E_{V}$ and
\item
[(b)] $x\cdot y=0$ if both $x,y\in V$ and $[x,y,z]=0$ if any two of $x,y,z$
lie in $V.$\\
\end{enumerate}

A particular instance where $V=\mathfrak{B}$ and we set $D(u,v)=D_{u,v},
\delta(u)=\delta_u$ the pair $(D,\delta)$ is a representation of Bol
algebra with companion $\Delta(u,v)=[u,-,v]$ called $regular representation$
of $\mathfrak{B}$.

\begin{expl}
Let $(\mathfrak{B},[-;-,-],\cdot)$ be the Bol algebra with basis $(e_1,e_2)$
over a field of complex numbers, were $[e_1,e_2,e_1]=e_{1}$, $[e_2,e_1,e_2]=e_{2}$ and
$e_1\cdot e_2=e_2$. We recall that $det(u,v)$ is the determinant of the
pair of vectors $(u,v)$ with $u=u_1e_1+u_2e_2$ and $u=u_1e_1+u_2e_2.$ Note that this Bol
algebra arise from the classification of two-dimensional Bol algebras obtained by Kuz'min and Zaidi in \cite{kuza}. We set

$ D(u,v)=\begin{pmatrix}
-det(u,v)& 0 \\
0& det(u,v)
\end{pmatrix}$

$ \delta(u)=\begin{pmatrix}
0 & 0 \\
u_2 & -u_1
\end{pmatrix}$

$ \Delta(u,v)=\begin{pmatrix}
-u_1v_2 & u_1v_1 \\
  u_2v_2  & -u_2v_1
\end{pmatrix}.$\\
It is clear that $(D,\delta,\Delta)$ is a regular representation of $\mathfrak{B}$.

\end{expl}

Now let $(\rho,\delta,\Delta)$ and $(\rho',\delta',\Delta')$ be two representations
of Bol algebra $\mathfrak{B}$ on $V,$ a morphism of the representation
$(\rho,\delta,\Delta)$ to a representation $(\rho',\delta',\Delta')$ is a linear
map $f:V\longrightarrow V$ such that $\rho'=f\rho$, $\delta'=f\delta$ and
$\Delta'=f\Delta.$ Clearly the composition of morphisms of representations is a
morphism of representations. The collection of all representations and their
morphisms forms a $K$-linear category
denoted by $Rep(\mathfrak{B})$ and called the category of representations of Bol
algebra $\mathfrak{B}$.\\
 We consider $Z_1(\mathfrak{B})=\bigcap\limits_{y\in \mathfrak{B}} ker(-\cdot y)$
and $Z_2(\mathfrak{B})=\bigcap\limits_{y,z\in \mathfrak{B}} ker[-;y,z]$, the center
of Bol algebra is $Z(\mathfrak{B})=Z_1(\mathfrak{B})\cap Z_2(\mathfrak{B}).$ It is
simple to see that, the kernel of the operation $<\rho,\delta>$ given by
$Ker<\rho,\delta>=\{x\in \mathfrak{B}/\rho(x,\mathfrak{B})+\delta(x)=0\}$ is the
center  of $\mathfrak{B}.$

\section{Engel's theorem for Bol algebras }
Before giving the Engel's theorem, we first need to define and characterize the
nilpotent representations.\\
A representation $(\rho,\delta,\Delta)$ of Bol algebra $\mathfrak{B}$ in $V$
is nilpotent if for all $x,y,z\in \mathfrak{B}$, $\rho(x,y),\delta(x)$  and $\Delta(x,y)$
are nilpotent endomorphisms;  that is if there is a positive integer $n$ such that
$(\rho,\delta,\Delta)^{n}=0.$
Let $(\rho,\delta,\Delta)$ be a representation of $\mathfrak{B}$ in $V,$ we define the triplet
$(ad_{\rho},ad_{\delta},ad_{\Delta})$ as follows: $ad_{\rho}(x,y)=[\rho(x,y),-]$,
$ad_{\delta}(x,y)=[\delta(x),-]$ and $ad_{\Delta}(x,y)=[\Delta(x,y),-].$

\begin{prop} With the above notations,
the pair $(ad_{\rho},ad_{\delta})$ is a representation of Bol algebra  $\mathfrak{B}$ in a
vector space $V$ with companion $ad_{\Delta}.$
\end{prop}

\begin{proof}
The objective is to show that $(R_1)$, $(R_2)$ and $(R_3)$ are satisfied.
Let $a,b,u,v\in \mathfrak{B}$ and $f\in EndV.$
We have
$$\begin{tabular}{lclll}
$[ad_{\rho}(a,b),ad_{\rho}(u,v)](f)$&=&$[ad_{\rho}(a,b),[ad_{\rho}(u,v),f]]$\\
&=&$[\rho(a,b),[\rho(u,v),f]]-[\rho(u,v),[\rho(a,b),f]]$\\
&=&$[[\rho(a,b),\rho(u,v)],f]$\\
&=&$[\rho(a,b),\rho(u,v)]f-f[\rho(a,b),\rho(u,v)]$\\
&=&$\rho([a;u,v],b)f+\rho(a,[b,u,v])f-f\rho([a;u,v],b)-f\rho(a,[b,u,v])$\\
&=&$(ad_{\rho}([a;u,v],b)+ad_{\rho}(a,[b,u,v]))(f)$
\end{tabular}$$

Then $(R_2)$ holds.
In other hand we have $[ad_{\rho}(a,b),ad_{\rho}(u,v)]=ad_{\rho}([a;u,v],b)+ad_{\rho}(a,[b,u,v])$

$$\begin{tabular}{lcl}
$ad_{\rho}(a,b)(f)$&=&$[{\rho}(a,b),f]$\\&=
&$\rho(a,b)f-f\rho(a,b)$\\
& =&$\Delta(a,b)f-\Delta(a,b)f-\Delta(b,a)f+f\Delta(b,a) $\\
&=& $[\Delta(a,b),f]-[\Delta(b,a),f]$\\
&= &$((ad_{\Delta}(a,b)-ad_{\Delta}(b,a))(f)$
\end{tabular}$$
Therefore we have the desire equality $ad_{\rho}(a,b)=ad_{\Delta}(a,b)-ad_{\Delta}(b,a).$
This shows that $(R_1)$ is satisfied.
Finally, we have for all $f\in EndV,$

$$\begin{tabular}{lcl}
$[ad_{\rho}(a,b),ad_{\delta}(u)](f)$&=&$ad_{\rho}(a,b)ad_{\delta}(u)f-ad
_{\delta}(u)ad_{\rho}(a,b)f$\\&=
&$[\rho(a,b),[\delta(u),f]]-[\delta(u),[\rho(a,b),f]]$\\
&=&$ [[\rho(a,b),\delta(u)],f]$\\
&=&$[\delta([u;a,b])+
\delta(a\cdot b)\delta(u)+ \Delta(a\cdot b,u),f]$\\
&=&$ [[\rho(a,b),\delta(u)],f]$\\
&=&$[\delta([u;a,b]),f]+
[\delta(a\cdot b)\delta(u),f]+[\Delta(a\cdot b,u),f]$\\
&= &$(ad_{\delta}([u;a,b])+
ad_{\delta}(a\cdot b)ad_{\delta}(u)+ad_{\Delta}(a\cdot b,u))(f)$
\end{tabular}$$
Thus $[ad_{\rho}(a,b),ad_{\delta}(u)]=ad_{\delta}([u;a,b])+
ad_{\delta}(a\cdot b)ad_{\delta}(u)+ad_{\Delta}(a\cdot b,u)$
and the desire conclusion follows, that is $(R_3)$ is verified.
\end{proof}
\begin{definition}
The representation $(ad_{\rho},ad_{\delta},ad_{\Delta})$ is
called the adjoint representation of $(\rho,\delta,\Delta).$
\end{definition}
Now we give the link between nilpotent representation and adjoint
representation. The above result arises to the representations of
Lie algebras.
\begin{lem}
Let $(\rho,\delta,\Delta)$ be a representation of Bol algebra
on the vector space $V.$ If $(\rho,\delta,\Delta)$ is nilpotent,
then its adjoint representation is also nilpotent.

\end{lem}
\begin{proof}
Let $(\rho,\delta,\Delta)$ be a nilpotent representation of Bol algebra,
and $(ad_{\rho},ad_{\delta},ad_{\Delta})$ its adjoint representation.
Then there exists a positive integer $p$ such that
$(\rho)^{p}=0$,$(\delta)^{p}=0$ and $(\Delta)^{p}=0$. If $\sigma$ is one
of the map $\rho,\delta$,or $\Delta$ it is clear that $ad_{\sigma}=
l_{\sigma}+h_{\sigma}$ where $l_{\sigma}$ and $h_{\sigma}$ are nilpotent.
we have $(ad_{\sigma})^{2p-1}=
(l_{\sigma}+h_{\sigma})^{2p-1}=0.$ Hence the result.

\end{proof}
Now we are in position to prove our first main theorem.

\begin{theo}
Let $\mathfrak{B}$ be a finite dimensional Bol algebra over a field $K$
and $\mathcal{R}$ consists of nilpotent representations of Bol algebra
$\mathfrak{B}$ in a finite dimensional space $V$. Then there exists a
vector $v\in V^{3}$, $v\neq 0$ such that $(\rho,\delta,\Delta)(v)= 0$ for
all $(\rho,\delta,\Delta)\in
\mathcal{R}.$

\end{theo}

\begin{proof}
We agree that
$(\rho,\delta,\Delta)(v)=(\rho(v_1),\delta(v_2),\Delta(_3))$, where
$v=(v_1,v_2,v_3),$ that is we identify $(\rho,\delta,\Delta)$ by $(\rho(a,b),\delta(a),\Delta(a,b))$
for all $a,b$ in $\mathfrak{B}$.
It is clear that $\mathcal{R}$ is a subspace of $(Env)^{3}$ and we can define
on it the following bracket $[(f,g,h),(f',g',h')]=([f,f'],[g,g'],[h,h']).$
$(\mathcal{R},[-,-])$ is a  Lie algebra.\\
The proof of the theorem goes by induction on $dim\mathcal{R}.$ When
$dim\mathcal{R}=1,$ since $\mathcal{R}$ is generated by a single
nilpotent representation then the claim is immediate.\\
Suppose now that the claim is true for all subalgebras of nilpotent
representations spaces of dimension less than $dim\mathcal{R}\geq 1.$\\
Since, $dim\mathcal{R}\geq 1,$ we have a proper Lie subalgebra
$L\subseteq \mathcal{R}.$ We can choose $L$ to be a maximal subalgebra.
We show before continuing that, $L$ has a codimension one in $\mathcal{R}$
 and $L$ is an ideal.\\ $L$ acts via the adjoint operator on $\mathcal{R}$
and $L$. In the latter case, since $dimL< dim\mathcal{R}$, we know by
Engel's theorem apply for $L$, that there
exists a nonzero element $\bar{r}\in \mathcal{R}/L$ such that $[l,\bar{r}]=0$
$\overline{(\rho,\delta,\Delta)}\in \mathcal{R}/L$ and
$[(l_1,l_2,l_3),\overline{(\rho,\delta,\Delta)}]=\bar{0}$ for $
(l_1,l_2,l_3)\in L.$ We know that $\overline{(\rho,\delta,\Delta)}=
(\rho,\delta,\Delta)+L;$ then $(\rho,\delta,\Delta)\in \mathcal{R}-L.$
It follows that $[K(\rho,\delta,\Delta)+L,L]\subseteq L.$
Moreover $[K(\rho,\delta,\Delta)+L,K(\rho,\delta,\Delta)+L]\subseteq L.$
These imply that $K(\rho,\delta,\Delta)+L$ is a Lie subalgebra of $\mathcal{R}$, and
contains $L$ as an ideal. By maximality  of $L$, it follows that $Kr+L=
\mathcal{R},$ so we are done.\\
Now we define the vector space $\mathbf{W}=\{w\in V^{3}/ Lw=0\}.$
Let $w=(w_1,w_2,w_3)\in \mathbf{W}$ and $(\rho,\delta,\Delta)\in L,$
then $(l_1,l_2,l_3)(\rho,\delta,\Delta)(w)=0$ for all $(l_1,l_2,l_3)\in L.$
Other we have

$$\begin{tabular}{lcl}
$(l_1,l_2,l_3)(\rho,\delta,\Delta)(w)$&=&$(\rho,\delta,\Delta)
(l_1,l_2,l_3)(w)+[(l_1,l_2,l_3),(\rho,\delta,\Delta)](w)$\\&=
&$[(l_1,l_2,l_3),(\rho,\delta,\Delta)](w)$
\end{tabular}$$
and $[(l_1,l_2,l_3),(\rho,\delta,\Delta)]\in L.$ Since $L$ is an ideal,
we have also $[(l_1,l_2,l_3),(\rho,\delta,\Delta)](w)=0.$\\
Now we have $\mathcal{R}=K(\rho,\delta,\Delta)+L$ for some
$(\rho,\delta,\Delta)\in L.$ We know that $(\rho,\delta,\Delta)$ is a
nilpotent operator on $\mathbf{W},$ so $ker(\rho,\delta,\Delta)\cap
\mathbf{W}\neq 0.$ Let $v=(v_1,v_2,v_3)\in ker(\rho,\delta,\Delta)\cap
\mathbf{W}$ such that $v\neq 0;$ then any element of $L$ and $r$
annihilates $v.$

\end{proof}

\section{An extension of Ado-Iwasawa theorem to Bol algebras}
Let $L$ be a finite-dimensional Lie algebra over a field $K.$
The classical Ado-Iwasawa theorem asserts the existence of a
finite-dimensional $L$-module which gives a faithful representation
of $L.$ However, Filippov proved in \cite{fil} showed that this
theorem does
not hold for Malcev algebras, that is homogeneous Bol algebras.
Thus it is not hold for general Bol algebras.\\
For the Lie algeras, the Poincar\'e-Birkhoff-Witt theorem says that
any Lie algebra $L$ is a subalgebra of $A^{-}$ for some unital
associative algebra $A.$ In the case that $L$ is finite dimensional,
the Ado-Iwasawa theorem says that $A$ can be taken finite dimensional
 too. This
extension of Ado-Iwasawa theorem was established for the Malcev
algebras in \cite{izsh}.
There is a version of the Poincar\'e-Birkhoff-Witt theorem for Bol
algebra proved in \cite{pi}. Now let
$\mathfrak{B}$ be a Bol algebra, it is shown in \cite{pi} that
there is an alternative algebra $A$ and an injective map
$\mathfrak{B}\longrightarrow RN_{alt}(A)$, where
$RN_{alt}(A)=\{a\in A/ (x,a,y)=-(x,y,a)\}$ is the generalized
right alternative nucleus. In this section
we prove that if $\mathfrak{B}$ is a finite-dimensional Bol algebra
then $A$ can be taken finite dimension too. Our second main result is
the following.

\begin{theo}
Let $\mathfrak{B}$ be a finite-dimensional right Bol algebra over a field
of characteristic $\neq 2,3$. Then there exists a unital finite-dimensional
algebra $A$ and a monomorphism of Bol algebra $j:\mathfrak{B}\longrightarrow
 RN_{alt}(A).$
\end{theo}

\begin{proof}
Let $\mathfrak{B}$ be a Bol algebra, according to P\'erez-Izquierdo in \cite{pi},
there exists a linear map $j:\mathfrak{B}\longrightarrow
RN_{alt}(U(\mathfrak{B})),$ $a\longmapsto a$  such that $j(a\cdot b)=ab-ba$ and
$j(a,b,c)=(ab)c-(ac)b-[b,c]a$, where $U(\mathfrak{B})$ is the universal
enveloping algebra of $\mathfrak{B}.$  Since $RN_{alt}(U(\mathfrak{B}))$
is closed under the binary product $[-,-]$ given by the commutators and the
ternary operation $[a,b,c]=(ab)c-(ac)b-[b,c]a$ for all $a,b,c$ in
$RN_{alt}(U(\mathfrak{B})).$ By the methods of \cite{pi},
$RN_{alt}(U(\mathfrak{B}))$ with the binary and ternary operations defined
above has the structure of Bol algebra. Thus $j$ is a monomorphism of Bol algebras.
Let $E_{\mathfrak{B}}$ be the Lie enveloping algebra of $\mathfrak{B}.$ Then
$E_{\mathfrak{B}}=E_{+}\oplus E_{-}$ is  the  $\mathbb{Z}_{2}$-gradation and
$E_{-}\approxeq \mathfrak{B}$ as vector space.
According to P\'erez-Izquierdo and Shestakov in \cite{izsh}, there exists a two
side ideal $\mathcal{I}\subseteq U(\mathfrak{B})$ of finite codimension.
Then $A=U(\mathfrak{B})/\mathcal{I}$ is a unital finite-dimensional algebra
and there exists an injective map $j:\mathfrak{B}\longrightarrow U(\mathfrak{B})$.
The injective map $j$ induces a monomorphism of  Bol algebras $j:\mathfrak{B}
\longrightarrow RN_{alt}(A).$

\end{proof}

\section{The category of representations of Bol algebra}
We give a relation between the category of representation of Bol algebra $\mathfrak{B}$ and the category of representations of its universal enveloping algebra. As immediate consequence, we show that the representation category of a Bol algebra is monoidal, or tensor category. We recall that the category of representations of Bol algebras is $Rep(\mathfrak{B})$, and the one of finite dimensional representations of  Bol algebra is $rep(\mathfrak{B})$. Let $A=(A,\cdot,\Delta,\epsilon)$ be a bialgebra, Mod($A$) means the category of left $A$-modules (ie., representations of $A$). If $U$, $V$ are left $A$-modules, then the tensor product becomes a left $A$-module with multiplication rule
$a\cdot (u\otimes v)=\Delta(a)\cdot (u\otimes v)$ for all $a\in A$, $u\in U$ and $v\in V$. The field $K$ is also a left $A$-module by $a\cdot \varsigma=\epsilon(a)\varsigma$. The category of left $A$-modules is equivalent to the category of $(A,A)$-bimodules. Any
$(A,A)$-bimodule can be considered as left module over $A\otimes A^{op}$, where $A^{op}$ is define on the same space as $A$, by new multiplication $x\cdot y=y\cdot x.$
We know in virtue of \cite{pi} that for a given Bol algebra $(\mathfrak{B},[-,-],[-,-,-])$ there exists a universal enveloping $U(\mathfrak{B})$ endowed with the structure of bialgebra, that is $(U(\mathfrak{B}),\cdot,\Delta,\epsilon)$ is a bialgebra. Analogously we denote $Rep(U(\mathfrak{B}))$ the category of representation of the bialgebra $(U(\mathfrak{B}),\cdot,\Delta,\epsilon).$
Now we state an equivalent characterization of the representation category $Rep(\mathfrak{B}).$ We prove our third main result.

\begin{theo}
The category of representations of Bol algebra $Rep(\mathfrak{B})$ is equivalent to the category of representations of its universal enveloping algebra $Rep(U(\mathfrak{B})).$
\end{theo}

\begin{proof}
We recall that $Rep(\mathfrak{B})$ is the category of modules over the Bol algebra $\mathfrak{B}.$ Following
the consideration of  Consuelo and Zelmanov see \cite{maze}, apply for the modules over Bol algebras, every
$\mathfrak{B}$-module has the form $E_{V}=\mathfrak{B}\bigoplus V$, where $V$ is a vector space over a field $K$
and $E_{V}$
possesses
the structure of Bol algebra such that:

\begin{enumerate} \item
[(a)] $\mathfrak{B}$ is a sub-Bol
algebra of $E_{V}$, \item
[(b)] $V$ is an ideal of Bol algebra $E_{V}$ and
\item
[(b)] $x\cdot y=0$ if both $x,y\in V$ and $[x,y,z]=0$ if any two of $x,y,z$
lie in $V.$\\
\end{enumerate}

We define the multiplication $U(\mathfrak{B})\times V\longrightarrow V$ by $a\cdot x=\epsilon(a)\cdot x$.
We consider the following mapping defined from $Rep(\mathfrak{B})$ to Mod$(U(\mathfrak{B}))$ define on the objets by $F(E_{V})=V.$
The map $F$ is naturally extended on the morphisms. If $U$ and $V$ are the images of $E_{U}$ and $E_{V}$ under $F,$ in virtue of \cite{pi}
there exits a map $\mu:\mathfrak{B}\longrightarrow U(\mathfrak{B})\otimes U(\mathfrak{B})$ with $\mu(a)=a\otimes 1+1\otimes a.$
This implies that $U\otimes V$ is a $U(\mathfrak{B})$-module.

    Conversely, let $V$ be a $U(\mathfrak{B})$-module, in virtue of P\'erez-Izquierdo, see \cite{pi} there exist an injective map
$\eta:\mathfrak{B}\longrightarrow U(\mathfrak{B})$. We define the multiplication $\mathfrak{B}\times V\longrightarrow V$ by $a\cdot x=\eta(a)\cdot x$. Then $V$ has the structure of module. We set now the mapping $G$ from Mod$(U(\mathfrak{B})$ to $Rep(\mathfrak{B})$ by $G(V)=E_{V}.$ It remains to define the image of $U\otimes V$. Let $E_{U}$ and $E_{V}$ be two modules over $\mathfrak{B}$, We set $E=\mathfrak{B}\oplus U\otimes V.$
We define the binary  operation by $[a,u\otimes v]_{\otimes}=[a,u]\otimes v$; $[a,u\otimes v]_{\otimes}=[a,u]\otimes v$ and a ternary by $[a,b,u\otimes v]_{\otimes}=[a,b,u]\otimes v$;  $[a,u\otimes v,b]_{\otimes}=[a,u,b]\otimes v$ and $[a,b,u\otimes v]_{\otimes}=[a,b,u]\otimes v$ for all $a$, $b$ in $\mathfrak{B}$, $u$ in $V$ and $v$ in $V.$ We assume also that the restrictions of $[-,-]_{\otimes}$ and $[-,-,-]_{\otimes}$ on $\mathfrak{B}$ correspond respectively to the binary and ternary operations of $\mathfrak{B}$; and
 $x\cdot y=0$ if both $x,y\in U\otimes V$ and $[x,y,z]=0$ if any two of $x,y,z$
lie in $U\otimes V.$ \newline
It remains to show that $(E,[-,-]_{\otimes},[-,-,-]_{\otimes})$ is a Bol algebra, that is the conditions $(i)$-$(iv)$ hold.
By the definition, the condition $(i)$ is satisfied. Now let $x,y,z,\alpha,\beta$ in $\mathfrak{B}$; $u$ in $U$ and $v$ in $V.$
We have

\begin{eqnarray*}[x;y,u\otimes v]+[u\otimes v;x,y]+[y;u\otimes v,x]&=& [x;y,u]\otimes v+[u;x,y]\otimes v+[y;z,u]\otimes v
\\&=& ([x;y,u]+[u;x,y]+[y;z,u])\otimes v
\\&=&0,\end{eqnarray*}

this shows that $(ii)$ is true.\newline

Now let us show that $(iii)$ holds. We have

\begin{eqnarray*}[[x; y,u\otimes v];\alpha,\beta]&=& [[x; y,u]\otimes v;\alpha,\beta]
\\&=&[[x; y,u];\alpha,\beta]\otimes v
\\&=&([[x;\alpha,\beta];y,u]+[x;[y;\alpha,\beta],u]+[x;y ,[u;\alpha,\beta]])\otimes v
\\&=& [[x;\alpha,\beta];y,u\otimes v]+[x;[y;\alpha,\beta],u\otimes v]+[x;y ,[u\otimes v;\alpha,\beta]]\end{eqnarray*}

One can show that the above equality holds for any $x,y,\alpha,\beta$ stands for $u\otimes v.$ That is $(iii)$ holds.

Finally, we have

\begin{eqnarray*}[[u\otimes v; y];\alpha,\beta]&=& [[u; y]\otimes v;\alpha,\beta]
\\&=&[[u; y];\alpha,\beta]\otimes v
\\&=&([u;\alpha,\beta]\cdot y+
[u,[y;\alpha,\beta]]+ [[\alpha,\beta]; u,y]+[[u,y],[\alpha, \beta]])\otimes v
.\end{eqnarray*}

Thus $[[u\otimes v; y];\alpha,\beta]=[u\otimes v;\alpha,\beta]\cdot y+
[u\otimes v,[y;\alpha,\beta]]+ [[\alpha,\beta]; u\otimes v,y]+[[u\otimes v,y],[\alpha, \beta]].$ One can show this
equality for any $y,\alpha,\beta$ stands for $u\otimes v.$ This completes the proof.

\end{proof}

\begin{definition}
A monoidal (tensor) category  $(\mathcal{C},\otimes,\mathbf{1},\alpha,\lambda)$ is a category $\mathcal{C}$
equipped with tensor functor $\otimes:\mathcal{C}\times\mathcal{C}\longrightarrow \mathcal{C}$, with a fix objet
$\mathbf{1}$ (called the unit of a tensor category),
$\alpha:\otimes \circ (\otimes \times Id)\longrightarrow \otimes \circ (Id \times \otimes)$,
$\lambda:\mathbf{1}\otimes-\longrightarrow Id$, $-\otimes \mathbf{1}\longrightarrow Id$ are natural isomorphisms
such that the associativity and unitary constraints hold, or equivalently the pentagon and the triangle diagrams
are commutative, for more details see (\cite{mac},\cite{eti},\cite{Kas}).
\end{definition}
We can now give a special characterization of the category of representations of Bol algebra as a consequence of the
above proposition.

\begin{cor}
Every category of representations of Bol algebras is a monoidal category.

\end{cor}

\begin{proof}
It was proved by Kassel \cite{Kas} that $(A,\cdot,\Delta,\epsilon)$  is bialgebra if and only if the category Mod($A$)
is monoidal category. In virtue of Theorem 5.0.6, the category of representations of Bol algebra is equivalent to the
category of representations of its enveloping algebra endowed with bialgebra structure. Hence the category $Rep(\mathfrak{B})$
is monoidal.

\end{proof}

More recently it was proved by Huang and Torecillas in \cite{hua}, that the path coalgebra $KQ$ of a given quiver $Q$ always admits
a bialgebra structure. So the monoidal category arising from this quiver bialgebra is the category of
representations of the bialgebra $KQ$. This leads to the following conjecture.

\begin{conj}
Find necessary and sufficient conditions for the existence of quiver $Q$ such that the monoidal category arising from
quiver bialgebra $KQ$ is the category of representations of a Bol algebra over algebraically closed field $K$.

\end{conj}

A monoidal category is said to be finite, if it is equivalent to the category of finite dimensional comodules over
 the finite dimensional coalgebra. Thus the category $rep(\mathfrak{B})$ of finite dimensional representations is finite
 monoidal category. This is a particular case of tensor categories of Etingof, Gelaki, Niksky and ostrik studied
 in \cite{eti}. The particular case where $Q$ is a quiver without loops and $2$-cyles should leads to strong relation
 between Bol algebras and $cluster$ $algebras$ of Fomin and Zelevinsky, see (\cite{fz},\cite{fz1}) for more details.
In the same vein, it has been shown in \cite{Shau}  that if $A$ is a finite dimensional bialgebra, then $A$ is Hopf algebra if and only
if the category of finitely generated $A$-modules is rigid, that is finitely generate modules admit dual objets. This allows us to the following conjecture.

\begin{conj}
Find necessary and sufficient conditions for a finite dimensional Bol algebra to have Hopf algebra as universal enveloping algebra.
\end{conj}

\textbf{Representations of free Bol Algebra $Bol[X]$ of finite dimension}

Let $X = \{x_{1},x_2,...,x_ n\}$, we construct the set of binary-ternary monomials $BT[X]$, and we assume that
$BT[X]$ is closed under
$[-,-]$ and $[-,-,-]$. Let $BT[X] = \{\sum\limits_{\substack{i=1}}^{n}{\alpha_{i}x_{i}}| \alpha_{i}\in \mathbb{F}\}$
 be the space spanned by $X.$
We define the multiplication by the following rules: if $f=\sum\limits_{\substack{i=1}}^{n}{\alpha_{i}x_{i}}$, $g=\sum\limits_{\substack{j=1}}^{n}{\beta_{j}x_{j}}$ and
$ h=\sum\limits_{\substack{k=1}}^{n}{\gamma_{k}x_{k}}$ in $BT[X],$ then $[f,g]=\sum\limits_{\substack{i,j=1}}^{n}{\alpha_{i}\beta_{j}[x_{i},x_{j}]}$, $[f,g,h]=\sum\limits_{\substack{i,j,k=1}}^{n}{\alpha_{i}\beta_{j}\gamma_{k}[x_{i},x_{j},x_{k}]}$. The free Bol
algebra $Bol[X]$ is the free
binary-ternary algebra $BT[X]$ satisfying the identities $(i)$-$(iv)$. The Bol types of degree $m$ are always
to construct a product of degree $m$ in $Bol[X]$. For general construction and more details of the free Bol
algebra $Bol[X]$, we refer to (\cite{hep},\cite{per}). In \cite{per} it has been shown that any multilinear
identity $f$ of degree $m$ can
be written as a linear combination of multilinear monomials. We denote the Bol  types of degree $m$ by $B_1$,
$B_2,$...,$B_{b(m)}$, that is $f=f_1+..+f_{b(m)}$, where $f_k$ is a linear combination of polynomial having Bol
type $k$. Therefore the author regards $f$ as an element of $b(m)$ copies of $\mathbb{F}S_{m}$, where
$\mathbb{F}S_{m}$ is group algebra of the group of permutation $S_m$. Applying the representation
$\Phi_{\sigma}:\mathbb{F}S_{m}\longrightarrow Md_{\sigma}(\mathbb{F})$, ($\sigma$ partition of $m$) of $S_m$ to
$f$ we obtain
the representation matrix of f in partition $\sigma$: $(\Phi_{\sigma}(f_1)|\Phi_{\sigma}(f_2)|...|\Phi_{\sigma}(f_{b(m)}))$.
\newline

Now let $V$ be finite dimensional space, dim($V$)=$s$ and $\mathfrak{B}$ is a Bol algebras of dimension $n$. Give a
representation $(\rho,\delta,\Delta)$ of $\mathfrak{B}$ over the space $V$ is equivalent to give the matrix
$(D(u,v)|\delta(u)|\Delta(u,v))$, where $D(u,v)$, $\Delta(u,v)$ are $s\times n$ matrices and $\delta(u)$ is also
a $s\times n$ matrix. Hence the block matrix $(D(u,v)|\delta(u)|\Delta(u,v))$ is a $(3n)\times s$ matrix.\newline

In the special case where $\mathfrak{B}=Bol[X]$, $K=\mathbb{F}$ and $V=\mathbb{F}S_{m}$, with Bol types $B_1$,
$B_2,$...,$B_{b(m)}$ the representation matrix $(\Phi_{\sigma}(f_1)|\Phi_{\sigma}(f_2)|...|\Phi_{\sigma}(f_{b(m)}))$ of
$f$ corresponds to the matrix $\delta_{f}$, that is the expression
$\delta(f)=(\delta(f_1)|\delta(f_2)|...|\delta(f_{b(m)}))$. At this specific case mentioned by Peresi \cite{per},
the representation of element $f$ is understood as a the representation of Bol algebra $Bol[X]$ given by the matrix
$(D(f,0)|\delta(f)|\Delta(f,0))$.

\bigskip

Actually we recall the classification theorem of Kuz'min and Zaidi for two-dimensional Bol algebras in \cite{kuza}
which states as follows.
\begin{theo}[Kuz'min-Zaidi]
Every Bol algebra $\mathfrak{B}$ of dimension two over $\mathbb{R}$ has a canonical basis $(e_1,e_2)$ in which its
multiplication table is one of the following:
 \begin{enumerate}
 \item[I.] $[e_1,e_2]=0$, $[e_2,e_1,e_2]=\varepsilon_{1}e_1,$ $[e_1,e_2,e_1]=\varepsilon_{2}e_2$, where $(\varepsilon_{1},\varepsilon_{2})=(0,0)$, $(-1,0)$,$(1,0)$,$(1,-1)$,$(1,1)$,$(-1,-1)$
  \item[II.]
  $[e_1,e_2]=e_2$, $[e_2,e_1,e_2]=\varepsilon e_1,$ $[e_1,e_2,e_1]=\beta e_2$, where $\varepsilon=0,-1,1;$
  $[e_2,e_1,e_2]=e_2,$ $[e_1,e_2,e_1]=e_1.$

\end{enumerate}

\end{theo}

Now we are in position to prove our classification result for regular representations of the two-dimensional Bol algebras.

\begin{theo}
Every regular representation of two-dimensional Bol algebra $\mathfrak{B}$ over $K$ is up to equivalence of matrices given by one of the following matrices:
 \begin{enumerate}

  \item[(1)]

$ R_1(u,v)=\begin{pmatrix}
0&\varepsilon_1 det(u,v)&0&0&u_2v_2\varepsilon_1 &-u_2v_1\varepsilon_1\\
-\varepsilon_2 det(u,v)&0&0&0&-u_1v_2\varepsilon_2 &u_1v_1\varepsilon_2
\end{pmatrix}$

 \bigskip

 \item[(2)]
$ R_2(u,v)=\begin{pmatrix}
0&\varepsilon det(u,v)&0&0&u_2v_2\varepsilon &-u_2v_1\varepsilon\\
-\beta det(u,v)&0&u_2&-v_1&-u_1v_2\beta &u_1v_1\beta
\end{pmatrix}$

 \bigskip

  \item[(3)]
$ R_3(u,v)=\begin{pmatrix}
-det(u,v)&0&0&0&u_1v_2 &u_1v_1\\
0&det(u,v)&u_2&-v_1&u_2v_2 &-u_2v_1
\end{pmatrix}$

\end{enumerate}

\end{theo}

\begin{proof}
In virtue of classification theorem of Kuz'min and Zaidi \cite{kuza}, every Bol algebra of dimension two
is of type $(I)$ or of type $(II)$ by using the items of their theorem. \newline

We suppose in the first case that our Bol algebra is of type $(I)$, that is $\mathfrak{B}$ has a canonical basis $(e_1,e_2)$ in which its multiplication table is given by
  $[e_1,e_2]=0$, $[e_2,e_1,e_2]=\varepsilon_{1}e_1,$ $[e_1,e_2,e_1]=\varepsilon_{2}e_2$, where \newline
  $(\varepsilon_{1},\varepsilon_{2})=(0,0)$, $(-1,0)$,$(1,0)$,$(1,-1)$,$(1,1)$,
  $(-1,-1).$ \newline
Let $u$ and $v$ be the two vectors of $\mathfrak{B}$, with $u=u_1e_1+u_2e_2$ and $u=u_1e_1+u_2e_2.$ We have $D(u,v)(e_1)=u_{1}v_{2}[e_1,e_1,e_2]+u_{2}v_{1}[e_1,e_2,e_1].$ Since $[e_1,e_1,e_2]=-[e_1,e_2,e_1]$, we have

\begin{eqnarray*}D(u,v)(e_1)&=& -u_{1}v_{2}[e_1,e_2,e_1]+u_{2}v_{1}[e_1,e_2,e_1]\\&=& (-u_{1}v_{2}+u_{2}v_{1})\varepsilon_{2}e_2
\\&=& -det(u,v)\varepsilon_{2}e_2,\end{eqnarray*}
 We have also
\begin{eqnarray*}D(u,v)(e_2)&=& u_{1}v_{2}[e_2,e_1,e_2]+u_{2}v_{1}[e_2,e_2,e_1]\\&=& (u_{1}v_{2}-u_{2}v_{1})\varepsilon_{1}e_1
\\&=& det(u,v)\varepsilon_{1}e_1.\end{eqnarray*}
 Thus
$ D(u,v)=\begin{pmatrix}
0&\varepsilon_1 det(u,v)\\
-\varepsilon_2 det(u,v)&0
\end{pmatrix}.$\newline
Now we compute the matrix of $\Delta(u,v)$ as follows. We have

\begin{eqnarray*}\Delta(u,v)(e_1)&=& u_{1}v_{2}[e_1,e_1,e_2]+u_{2}v_{2}[e_2,e_1,e_2]\\&=& u_{2}v_{2}\varepsilon_{1}e_1-u_{1}v_{2}\varepsilon_{2}e_2
,\end{eqnarray*}

and
\begin{eqnarray*}\Delta(u,v)(e_2)&=& u_{1}v_{1}[e_1,e_2,e_1]+u_{2}v_{1}[e_2,e_2,e_1]\\&=& -u_{2}v_{1}\varepsilon_{1}e_1+u_{1}v_{1}\varepsilon_{2}e_2
,\end{eqnarray*}
hence
$ \Delta(u,v)=\begin{pmatrix}
u_{2}v_{2}\varepsilon_{1} & -u_{2}v_{1}\varepsilon_{1}\\
-u_{1}v_{2}\varepsilon_{2}& u_{1}v_{1}\varepsilon_{2}
\end{pmatrix}.$ Because $[e_1,e_2]=0,$ we have $\delta(u)=0.$ Therefore the bloc matrix $(D(u,v)|\delta(u)|\Delta(u,v))$
corresponds to the matrix $R_1(u,v).$

 The second case corresponds to Bol algebra  of type $(I)$, that is $\mathfrak{B}$ has a canonical basis $(e_1,e_2)$ in
 which its multiplication table is given by
  $[e_1,e_2]=e_2$, $[e_2,e_1,e_2]=\varepsilon e_1,$ $[e_1,e_2,e_1]=\beta e_2$, where $\varepsilon=0,-1,1;$
  $[e_2,e_1,e_2]=e_2,$ $[e_1,e_2,e_1]=e_1.$ \newline

If $[e_1,e_2]=e_2$, $[e_2,e_1,e_2]=\varepsilon e_1,$ $[e_1,e_2,e_1]=\beta e_2$, where $\varepsilon=0,-1,1;$ we use the
analogous methods as at the first case to get  $ D(u,v)=\begin{pmatrix}
0&\varepsilon det(u,v)\\
-\beta det(u,v)&0,
\end{pmatrix}$
  $ \delta(u)=\begin{pmatrix}
0  & 0 \\
u_2 & -u_1
\end{pmatrix}$
and
$ \Delta(u,v)=\begin{pmatrix}
u_{2}v_{2}\varepsilon & -u_{2}v_{1}\varepsilon\\
-u_{1}v_{2}\beta & u_{1}v_{1}\beta
\end{pmatrix}.$ Hence the bloc matrix $(D(u,v)|\delta(u)|\Delta(u,v))$
corresponds to the matrix $R_2(u,v).$

Finally, for$[e_1,e_2]=e_2$ and $[e_2,e_1,e_2]=e_2,$ $[e_1,e_2,e_1]=e_1,$ we have

\bigskip

$(D(u,v)|\delta(u)|\Delta(u,v))=\begin{pmatrix}
-det(u,v)&0&0&0&u_1v_2 &u_1v_1\\
0&det(u,v)&u_2&-v_1&u_2v_2 &-u_2v_1
\end{pmatrix},$ this end the proof.
\end{proof}

\textbf{Acknowledgement}

We aim at sending our sincere gratitude to the reviewers for their reports
on our paper and their availability. The second author thanks the IHES for
hospitality through the Launsbery Foundation and, the IMU for travel
support throughout a grant from the Simons Foundation during the writing
of this paper.

\end{document}